\documentclass[
aip, 
jmp,
amsmath, amssymb, 
reprint,
]{article}

\usepackage{graphicx}
\usepackage{makeidx}
\usepackage{multirow}
\usepackage{amssymb}
\usepackage{amsfonts}
\usepackage{amsmath}
\usepackage{amsthm}
\usepackage{natbib}

\usepackage{graphicx}
\usepackage{dcolumn}
\usepackage{bm}

\usepackage[utf8]{inputenc}
\usepackage[T1]{fontenc}
\usepackage{mathptmx}
\usepackage{amsfonts}
\usepackage{bbm}
\usepackage[bottom]{footmisc}

\begin{document}

\title{On the diversity and similarity of mathematical models in science}

\author{Inge S. Helland
\\ Department of Mathematics, University of Oslo \\ P.O.Box 1053 Blindern, N-0316 Oslo, Norway
\\ Orcid: 0000-0002-7136-873X\\ Telephone: (047)93688918\\ E-mail: ingeh@math.uio.no}

\date{\today}

\maketitle

\begin{abstract}
In this article, the notion of a mathematical model in science is attempted to be enlightened from several points of view. In particular, it is shown that mathematical models are introduced differently and used differently in different areas of science. In the present article the use of models in statistics is taken as a basis, but links are given to several other areas. Two particular such links are described in some detail here: First the link connected to the chemometricians' Partial Least Squares algorithm, a link that now has been generalized to the more sophisticated envelope model. Secondly, statistics, as is well known, relies heavily on making decisions, and it may also in certain cases be connected to a process of taking measurements. A mathematical model for these two activities, connected to the mind of an observer, is introduced. This model is taken further, and shown to be important in a new proposal for a foundation of quantum theory. Quantum mechanics, as seen from this point of view, is described in some detail. The discussion here is close to the discussion in a recent book and in a recently published article by the author.

\end{abstract}

\bigskip

\section{Introduction}
\label{Sec1}

Biometrics, chemometrics, psychometrics and econometrics are sciences that are closely connected to statistics, but are more focused on specific applications. Communication across scientific borders is important, and to some extent such a communication requires a common language.

The language of mathematics is common to all natural sciences. The present article will focus on this language, and on the use of mathematical models in science. Even though there are cultural barriers, my thesis is that one should be able to formulate mathematical models in such a way that they - at least after some effort -  could be understood by all serious researchers, in all disciplines.
 The word `understand' is quite demanding. My late colleague Emil Spj\o tvold once remarked that things can be understood on many different levels. This applies also in this case.
 
 Take one important example, quantum mechanics. In textbooks and in research papers, quantum theory is formulated in a very abstract language, and using this language, physicists are able to compute, but the issue of interpretation is very controversial, and researchers from other disciplines have great difficulties of understanding what is going on.
 
 In fact, quantum mechanics was originally formulated in two different languages, the wave mechanics of Schr\"{o}dinger and the matrix mechanics of Heisenberg. By mathematical abstraction these two were unified, and that resulted in a mathematically clean, but conceptual difficult theory theory, which has been extremely difficult to interpret. The relevant Wikipedia article contains more than 16, partially mutually excluding, interpretations.
 
In the book Helland (2021) and the article Helland (2022a), a new basis is proposed, a basis which is hopefully easier to understand for researchers outside the quantum community. As I see it, we can take as our point of departure a set of conceptual variables $\theta$ in the mind(s) of the investigator(s), and these variables are introduced through questions of the form `What will $\theta$ be if we measure it?'. In the discrete case, which is the focus of many textbooks and research papers, one can obtain sharp answers of the form `$\theta = u$'. For certain well-defined conceptual variables, these questions and answers together define the relevant quantum states. This should be compared to the ordinary quantum language, where a state is defined as a unit vector with arbitrary phase in an abstract, complex Hilbert space.

Below, I will discuss the concept of a mathematical model in general. I will cover several disciplines, but the story starts with the interplay between chemometrics and statistics. Since so many different areas are covered, I will attempt to be brief on each single discipline. This may imply some simplifications. In particular, the choice of references here is somewhat subjective.

The important set of deterministic models will not be covered here. Nor will I mention stochastic models in, say, population dynamics and genetics. I will concentrate on models that, in some way of other, may be associated with areas connected to statistical inference.

The plan of the paper is as follows: In Section \ref{Sec2} I give a general discussion of cultures within science. Then, as a starting point, mathematical models used in statistical inference is described in Section \ref{Sec3}. The Partial Least Squares `model' as developed in Chemometrics in the 1980's is discussed in Section \ref{Sec4}, while Section \ref{Sec5} gives a link between Section \ref{Sec3} and Section  \ref{Sec4}. Dennis Cook's envelope model, which can be closely associated with this, is considered in Section \ref{Sec6}. In Section \ref{Sec7} a brief discussion of the connection between statistical models and machine learning is given, while Section \ref{Sec8} contrasts briefly statistical models and models connected to causality. Decisions and measurements are important both in connection to statistics and in connection to the foundation of quantum theory, and models for these are introduced in Section \ref{Sec9}. Then in Section \ref{Sec10} the elements of a foundation of quantum mechanics, taking the model notion as a basis, is described. Section \ref{Sec11} gives some general discussion points.

\section{On science and cultures within science}
\label{Sec2}

My definite starting point is that science is important in society. There are extremely many problems facing the world now: climate, health problems, poverty problems, refugee problems, international conflicts and the existence of very dangerous weapons. For all these problems one should try to device rational solutions that at the same time satisfy good ethical standards. In particular, rational and good ethical decisions should be made by national and international leaders. Unfortunately, in the process of making decisions, we may all be limited. (In a concrete physical situation where communication is important, this has been argued for in Helland, 2022b.) So ideally, to arrive at good decisions, in addition to other factors, our political leaders should attempt to get good, rational inputs. Ideally again, such inputs could in concrete cases be given by a joint effort by groups of scientists.

In many cases this requires time, and it may also meet practical difficulties. But it is important that scientists are focused in such cases, and have the necessary tools for obtaining the relevant scientific knowledge.

In choosing tools, a scientist must make choices that sometimes are uncertain. It is then important that he can get inspiration from colleagues. Such set of colleagues may then be a part of a larger community. What we do not often think about, is that such communities form a joint culture, and that elements of such a culture may have arbitrary or historical elements. Sometimes a scientist's choice of tools may depend on what inputs he has obtained from his culture.

An example can be the choice of regression method in the case of collinear data. Statisticians have many tools here, variable selection, ridge regression, principal component regression and so on. In addition, chemometricians have developed their own tool, partial least squares (PLS), which empirically seems to to have good properties, and has now been applied in a great variety of fields (Mehmood and Ahmed, 2015).

What is culture? According to the author and philosopher Ralph D. Stacey it is a set of attitudes, opinions and convictions that a group of people share, about how one should act towards each other, how things should be evaluated and done, which questions that are important and which answers that may be accepted. The most important elements in a culture are unconscious, and cannot be forced upon from the outside.

\section{Models, probabilities and procedures in statistics}
\label{Sec3}

Most of what we do in statistics is based on probability models for the data, but there are exceptions (Breiman, 2001). I will disregard nonparametric models here, and assume continuos data, so that the model in general is given by a probability density
\begin{equation}
f_\theta (x),
\label{1}
\end{equation}
where $x$ is the data, and $\theta$ is the full parameter. It is important that the model is selected by the researcher. The parameter $\theta$ exists in some fundamental way in the researcher's mind, and is related to the research question that he wants to investigate. More precisely, the research question is often posed in terms of a subparameter $\eta =\eta (\theta)$.

It is well known that statistical methods may be classified into frequentist methods, Bayesian methods, and - as a perhaps unusual way of thinking - fiducial methods. To some extent this divides statisticians into different schools, but firstly, there are connections between the different set of methods, and secondly, one single person may also use methods from different schools.

The interpretation and use of the probability concept varies somewhat between different schools. The Bayesians have their prior distributions, and, from this and the statistical model, posterior distributions. A frequentist will sometimes distinguish between aleatoric probabilities - probabilities of data, given some assumed properties of the world - and epistemic probabilities - probabilities of some properties of the world, perhaps given some data. I will come back to the latter below. The former are the bases for constructing statistical models.

From a Bayesian perspective, everything can be modeled by using ordinary probability models. The priors and the posteriors both follow Kolmogorov's axioms and all other probability rules that we are used to. In particular, the posterior distributions are epistemic distributions that can be seen as ordinary probability distributions.

In Schweder and Hjort (2016), confidence distributions were proposed from a frequentist point of view as a method of making `probability' statements about aspects of the world, given some data from an experiment. Here is the definition from Schweder and Hjort (2016):
\bigskip

\textbf{Definition 1} \textit{A nondecreasing right-continuous function of the one-dimensional $\eta$, depending on the data $Y$, say $C(\eta ,Y)$, is the cumulative distribution function for a confidence distribution for $\eta$ provided that it has a uniform distribution under the statistical model, whatever the true value of the full parameter $\theta$.}
\bigskip

More concretely, this can be obtained as follows: Take a class of one-sided confidence intervals $(-\infty , \eta_\alpha ]$ corresponding to the confidence coefficient $\alpha$. Then, we can find $C$ from $\eta_\alpha =\eta_\alpha (Y)= C^{-1}(\alpha)$, so that $C$ also depends on the data $Y$:
\begin{equation}
P_\theta [\eta\le \eta_\alpha]=P_\theta[C(\eta )\le C(\eta_\alpha )]=P_\theta [C(\eta)\le\alpha]=\alpha .
\label{2}
\end{equation}

For a scalar continuous parameter, Fisher's fiducial distribution is essentially the same as the confidence distribution. More generally, a fiducial distribution can be found (Taraldsen and Lindqvist, 2021) from a data generating model $y=y(\eta ,u)$, where $\eta=\eta(\theta )$ is a parameter, and $u$ has a probability distribution depending upon $\theta$. The model is called simple if it has a unique solution $\eta =\eta(y,u)$. The fiducial distribution is then the distribution of $\eta$ for fixed $y$.

Fisher's difficulties began when he tried to generalize his concept to the multivariate case. In the last few years, the concept of a fiducial theory has been taken up again, see references in Taraldsen and Lindqvist (2021). Older results by Fraser (1964) has tied the fiducial distribution to a Bayesian distribution with a group-invariant prior.

A concept due to Fisher (1973) is the likelihood function $L$ the probability density (\ref{1}), seen as a function of $\theta$. Maximum likelihood estimation is used in many applied fields. Less well known is the likelihood principle, saying that observations with the same, or proportional, likelihood always give the same experimental evidence. In Helland (2021) I advocate a version where the experimental context is fixed, which in my opinion makes the principle less controversial.

A well known statistical model with a large number of applications is the regression model
\begin{equation}
y=\mu + X\beta+e,
\label{3}
\end{equation}
where $y$ is an $n$-dimentional vector of response variables, $X$ is an $n\times p$ matrix of explanatory variables, and $e$ is a random error term. When there is no collinearity, the maximum likelihood estimator under normality is
\begin{equation}
\widehat{\beta}=(X'X)^{-1}X'y .
\label{4}
\end{equation}

When there is collinearity, in particular, if $p>n$, many different estimators have been proposed by statisticians. In general, our purpose might be prediction of a new $y$-value $y_0$ from a set of $x$-values $x_0$:
\begin{equation}
\widehat{y_0}=\widehat{\mu}+ \widehat{\beta}'x_0 .
\label{5}
\end{equation}

\section{A certain `model' from chemometrics}
\label{Sec4}

During the 1980's chemometricians developed an algorithm which was intended to solve the above regression problem in cases with collinearity: The partial least squares, or PLS algorithm. The algorithm used had its origin in Herman Wold's general system analysis (J\"{o}reskog and Wold, 1982). Briefly, a theoretical version of it goes as follows, where the starting point is an $x$-vector $x$ (the column vector corresponding to a typical row of $X$ in (\ref{3}), and a scalar $y$ (the corresponding component of the vector $y$ in (\ref{3}):
\bigskip

(i) Define starting values for the $x$ residuals ($e_a$) and $y$ residuals ($f_a$):
\[e_0 =x-\mu_x , \  f_0 = y-\mu_y \]

For $a=1,2,...$, do steps (ii)-(iv) below:

(ii) Introduce scores (t) as linear combinations of the $x$ residuals, using covariates with $y$ as weights ($w$):
\[t_a =e_{a-1}'w_a ,\  w_a = \mathrm{Cov} (e_{a-1}, f_{a-1}).\]

(iii) Determine $x$ loadings ($p_a$) and $y$ loadings ($q_a$) by least squares:
\[p_a =\mathrm{Cov} (e_{a-1}, t_a )/\mathrm{Var}(t_a),\  q_a =\mathrm{Cov} (f_{a-1}, t_a )/\mathrm{Var}(t_a).\]

(iv) Find new residuals:
\[e_a =e_{a-1}-p_a t_a ,\  f_a =f_{a-1} - q_a t_a .\]

This whole procedure, including the construction of scores, loadings and weights, was called `soft modelling'. At each step $r$ a bilinear representation resulted from this:
\[x=\mu_x+p_1 t_1+...+p_r t_r +e_r,\ y=\mu_y +q_1 t_1+...+q_r t_r+f_r .\]

The PLS \emph{algorithm} is then used for empirical prediction, a method that has become very popular in many applied fields (Mehmood and Ahmed, 2015):

Use the data $(X,y)$, where $y$ now is a vector as in (\ref{3}). From these data estimate expectations, variances and covariances; otherwise use the procedure to find estimates of scores, loadings and weights. Stop the algorithm at the step $r=m$, found by cross-validation. Predict new $y$-values by:
\[\widehat{y_0} = \bar{y}+\widehat{q_1}\widehat{t_1}+...+\widehat{q_m}\widehat{t_m}.\]

All this was met with some sceptisism by some statistians; other statisticians chose to neglect the whole development. For those who were engaged, a crucial debate resulted from this; see for instance Martens (1987) and Schweder (1987).

Essentially two problems were the reasons for the statistical sceptisism. Firstly, the socalled soft modelling was not related to the `proper' modelling that statisticians were used to. A model should distinguish sharply between data on the one hand, and parameters on the other hand. Secondly, the algorithm feature and the strongly non-linear characteristic of the procedure, made it very difficult to develop ordinary statistical properties of the method, so that it could be compared to other methods.

The first problem was essentially solved in 1990; the second problem has only begun to reach a solution in the last few years; see below.

\section{The link}
\label{Sec5}

In Helland (1990) a proper statistical PLS model was proposed. Take as a point of departure the model (\ref{3}), where the $x$-variables are assumed to be random variables, for simplicity centered to zero expectation. Furthermore, we assume that the $x$-variables have a common covariance matrix $\Sigma_x$ and are uncorrelated with the errors $e$. 

The crucial parameters are then $\Sigma_x$ and $\beta$. It was shown in Helland (1990) that the scores, loadings and weights in the empirical PLS algorithm could be seen as estimates of scores, loadings and weights in a corresponding population algorithm, made to stop automatically at step $m$ ($e_m =f_m =0$), $(m< \mathrm{min} (n,p))$. Furthermore, this automatic stop-restriction could be explicitly formulated in two equivalent ways in terms of $\Sigma_x$ and $\beta$:

(1) Consider the expansion $\beta =\sum_{i=1}^p \gamma_i d_i$, where $d_1 ,...,d_p$ is a complete set of eigenvectors of $\Sigma_x$. The PLSR model with $m$ factors is defined as the model where this expansion is reduced to exactly $m$ terms. One does not say anything about \emph{which} eigenvectors are involved in the expansion.

(2) Let $\sigma= \Sigma_x^{-1}\beta$, and consider the Krylov sequence $\sigma , \Sigma_x \sigma , \Sigma_x^2 \sigma ,\Sigma_x^3 \sigma,...$. The PLSR model with $m$ factors is the model where the space spanned by this sequence has dimension $m$.

These ideas were developed further in a few papers: In Helland (1992) the maximum likelihood estimator of $\bm{\beta}$ under the PLSR model was derived. In N\ae s and Helland (1993) the concept of relevant components, the eigenvectors correlated with $\beta$ (cf. formulation (1) above) was discussed in detail, and this concept was related to other regression methods in Helland and Alm\o y (1994). A group theoretical approach towards the PLSR model was discussed in Helland et al. (2012).

\section{Envelope models and PLS}
\label{Sec6}

From 2010 on, a new development took place. In the article Cook et al. (2010), the important general \emph{envelope model} was proposed and discussed, and in Cook et al. (2013) it was shown that a PLSR model - extended to the case of multivariate $y$ - was an important case of the envelope model. For the multivariate case, the regression coefficient $\beta$ in (\ref{3}) is replaced by a matrix $B$. The envelope model is - quite similarly to the PLS model, a restiction of the joint parameter $(\Sigma_x ,B)$. This can be formulated in several ways; one is the following: Assume that there exist an orthogonal matrix $(\Phi ,\Phi_0 )$, where $\Phi$ is $p\times m$, $\Phi_0$ is $p\times (p-m)$ and matrices $\eta, \Delta$ and $\Delta_0$ such that $B=\Phi\eta$ and
\begin{equation}
\Sigma_X =\Phi \Delta\Phi'+\Phi_0 \Delta_0 \Phi_0'.
\label{6}
\end{equation}
The first term in the expansion for $\Sigma_X$ corresponds to the covariance function for the `relevant' linear combinations of $x$-variables, the second term to the `irrelevant' linear combinations.

A large number of papers on the envelope model have been published by Dennis Cook and his collaborators in the last decade, and the theory is summarised in the book Cook (2018). As estimation in this model is concerned, the development is dominated by maximum likelihood estimation, which does not work for $p>n$. A Bayesian approach for envelope models has recently been discussed in Khare et al. (2017).

What is the situation now?  From a statistical point of view, much is still unknown about the PLS algorithm, but see the recent paper by Cook and Forzani (2019) and references there. In the latter paper, asymptotic expansions also in the case $p\rightarrow\infty$ were studied, and the algorithm was shown to have good statistical properties in many cases when many $x$-variables have non-zero regression coefficient (no sparsity).

Both for PLS and for the developments initiated by Dennis Cook and his students, it is important that ordinary statistical models can be said to lie behind the predictions done, and the more general applications. A statistical model depends on parameters, and these parameters must be said to have an existence either in the mind of a single researcher, or in the joint minds of a group of researchers. This is a notion that is the background for all statistical modeling, and a notion that we should hold on to, also when trying to approach other other areas of (natural) science.

\section{Machine learning models}
\label{Sec7}

Several statisticians have changed their focus in the last few years. Artificial intelligence, and particularly, machine learning has turned out to be very important in many applications. There are many books on machine learning. A good modern book on the computer approach as seen from statisticians' point of view, is Efron and Hastie (2021).

There is no reason to review the content of the book Efron and Hastie (2021) here, except to emphasize that it is firmly based on statistical models. Chapter 5 discusses models in general, Chapters 8 and 9 uses statistical models as a basis for the authors' treatment of early computer-age methods, while Chapter 16 on sparse modeling is contained in the Section `Twenty-First-Century Topics'. But the book also contains a Chapter (18) on  neural networks and deep learning.

An important paper on learning rules seen from a Bayesian perspective is Khan and Rue (2021). For those who are interested in the cultural aspect of modern statistical learning, in particular the merging of the model culture and the more algorithmic approach, I mention the brief article Bhadra et al. (2021). 

\section{Causality and statistical models}
\label{Sec8}

The book Pearl (2009) provides a systematic account of cause-effect relationships among variables or events. Modeling causality should be sharply distinguished from statistical modeling. Important notions in Pearl (2009) are directed acyclic graphs, Bayesian networks and the introduction of a do-statement.

But this distinction can be formulated in terms of probabilistic modeling, as discussed by Cox and Wermuth (2004). Take as a point of departure three sets of variables $B$ (background), $C$ (causes) and $R$ (response). These can be modeled in recursive form as
\begin{equation}
f_{RCB}=f_{R|CB}f_{B|C}f_B .
\label{7}
\end{equation}

A corresponding full graph of the three sets of variables will contain an arrow from  $B$ to $R$, an arrow from $B$ to $C$, and an arrow from $C$ to $R$. With this full graph, the distribution of $R$, given $C$ can be formed by marginalizing over $B$:
\begin{equation}
f_{R|C} = \int f_{R|CB} f_{B|C} db ,
\label{8}
\end{equation}
where $f_{B|C}=f_{BC}/f_C$. This is a \emph{statistical} way of modeling, taking full account of the notion of conditioning.

By contrast, in a \emph{causal} model the arrow from $B$ to $C$ is assumed to be missing. Then we have a case of intervention from $C$ to $R$, which can be written as
\begin{equation}
f_{R||C} = \int f_{R|CB}f_B db .
\label{9}
\end{equation}

The distinction between the two probability distributions $f_{R|C}$ and $f_{R||C}$ is crucial for any discussion of inference versus intervention, that is, statistics versus causality.

\section{Models for decisions and measurements}
\label{Sec9}

In the whole discussion above, it is assumed that the background for modeling, and in particular the set of statistical parameters, is something that belongs to the mind of a single researcher or to the joint minds of a group of communicating researchers. Let us first assume in general that the mind of some person $A$ contains variables $\theta, \xi, \lambda,...$. In Helland (2021) these are called conceptual variables.

Some of the conceptual variables may be decision variables: Let us assume that $A$ is faced with some concrete decision, having a finite number of prospects $\pi_1 ,\pi_2 ,...,\pi_r $. Then he can define for himself a decision variable, taking the values $1, 2,...,r$ such that $\xi =k$ is identified with the choice of the future prospect $\pi_k$ $(k=1,2,...,r)$. 

One situation lying behind much of the discussion in Helland (2021) is that of a measuring process. $A$ may be planning to do some measurement of a physical variable, say, either $\theta$ or $\eta$ or .... These variables may be assume to be connected to some concrete physical system, but a crucial assumption in Helland (2021) is that they, before, during and after a measurement, also belong to the mind of some person $A$. This is in agreement with the philosophy of Convivial Solipsism, as recently was proposed by Zwirn (2016) as part of a background for quantum theory: All aspects of the world must be seen from the point of view of the mind of some person. But people can communicate.

Think generally. $A$ is in some specific context, he has his physical environment, his own history, his plans about doing something, and so on. In this context he has several conceptual variables in his mind at some given point of time. Some of these are what I call accessible: It is possible by a measurement or in other ways to find values for them in some future. (Note that this also is meaningful for decision variables.)

We can define a partial ordering both among all conceptual variables and among the accessible ones: Say that $\theta$ is `less than or equal to' $\lambda$ if $\theta = f(\lambda)$ for some function $f$. I will assume that, if $\lambda$ is accessible, then $\theta=f(\lambda)$ also is accessible.

In many situations $A$ will have what I call maximally accessible conceptual variables. First look at a measurement situation. Then this may be achieved by variants of Heisenberg's inequality. In particular then, position and momentum for a particle at some time will both be maximally accessible. These variables are then also called, after Niels Bohr, complementary. In the last chapter of Helland (2021), the complementerity notion is attempted generalized to many different areas, also macroscopic.

But go back to the person $A$, and attempt to look at all his conceptual variables at time $t$ from a psychological point of view. We can make a model for this by assuming that there is an inaccessible variable $\phi$ such that all or a set of accessible variables can be seen as functions of $\phi$. In some way $\phi$ must be said to lie in the subconscious part of $A$'s mind, and it can not be determined by $A$ himself. A person $B$, knowing $A$ well and having some insight in practical psychology, may perhaps in some situations come up with a rough estimate of $\phi$.

As any model, this may be a simplification, but it turns out to be a fruitful simplification.

Relative to this model, and using Zorn's Lemma (in the partial ordering $\phi$ is a dominating conceptual variable), there certainly exist maximally accessible conceptual variables. This is what we need for the next Section.

The paper Yukalov and Sornette (2014) and several related papers lie the foundation of what is called Quantum Decision Theory, a quantum model for decisions made by a single person. In Yukalov et al. (2018, 2022) this theory is extended to joint decisions made by a group of persons. For this, exchange of information is crucial. In Yukalov et al. (2022), collective information shared by the whole group, is also taken into consideration. The decision process for a group of people may take time.

Yukalov and Sornette (2014) derive the quantum probabilities $p(\pi_i)$, where $\pi_i$ is the prospect, from what they call an utility factor $f(\pi_i )$ and an attraction factor $q(\pi_i)$: $p(\pi_i)=f(\pi_i)+q(\pi_i)$. Here $f(\pi_i)$ corresponds to the conscious part, modeled through an utility function and a coefficient expressing the actor's belief or confidence. The factor $q(\pi_i)$ corresponds to the unconscious part, and satisfies $-1\le q(\pi_i)\le 1$ and $\sum_i q(\pi) =0$. In Yukalov et al. (2018, 2022) there is in addition a factor $h(\pi_i)$ depending on the (joint) information.

Alternatively, quantum probabilities can be derived from Born's formula.

\section{Models as a fundation of quantum mechanics}
\label{Sec10}

Several authors have recently tried to rederive the axioms of quantum mechanics, which lead to an abstract theory, from simpler assumptions. As is well known, a basic axiom of quantum mechanics is that a physical system at some fixed time can be associated by a Hilbert space, a complex vector space with a scalar product that is complete in the norm derived from this scalar product, and that the pure states of this system can be defined as unit vectors with arbitrary phase in this Hilbert space. This works well as a basis for quantum calculation, but is very abstract, is difficalt to understand for outsiders, and has lead to many, partially conflicting, intepretations of the theory.

In Helland (2021, 2022a) this problem is approached by taking conceptual variables as the basic notion. The proofs needed for developing the theory here are somewhat complicated, using group theory and group representation theory, but from this approach, the conceptual basis  for the theory can be formulated in a relatively simple way.

The proof of these results relies on some technical theorems. The first of these can be formulated as follows:
\bigskip

\textbf{Theorem 1}
\textit{Consider a situation where there are two maximally accessible conceptual variables $\theta$ and $\xi$ in the mind of an actor or in the minds of a communicating group of actors. Make the following assumptions:}

\textit{(i) On one of these variables, $\theta$, there can be defined transitive group actions  $G$ with a trivial isotropy group and with a left-invariant measure $\rho$ on the space $\Omega_\theta$.}

\textit{(ii) There exists a unitary irreducible representation $U(\cdot)$ of the group behind the group actions $G$ defined on $\theta$ such that the coherent states $U(g)|\theta_0\rangle$ are in one-to-one correspondence with the values of $g$ and hence with the values of $\theta$.}

\textit{(iii) The two maximally accessible variables $\theta$ and $\xi$ can both be seen as functions of an underlying inaccessible variable $\phi\in\Omega_\phi$. There is a transformation $k$ acting on $\Omega_\phi$ such that $\xi(\phi)=\theta(k\phi)$.}

\textit{Then there exists a Hilbert space $\mathcal{H}$ connected to the situation, and to every accessible conceptual variable there can be associated a symmetric operator on $\mathcal{H}$.}
\bigskip

For those that are unfamiliar with the ordinary textbook-formulation of quantum theory, it must be repeated that one of the postulates stated in the textbooks is just that there in every physical situation exists a Hilbert space, and furthermore that every observable variable is associated with a socalled self-adjoint operator in this Hilbert space. Now self-adjoint is in most cases essentially the same as symmetric (Hall, 2013).

The point of the theorem is to derive this conceptual apparatus from simpler assumptions. First look at a measurement situation, as exemplified by measuring either position or momentum of a particle. As discussed in the previous section, these variables can also be seen as conceptual variables, existing in the mind of some observer. For simplicity, look at the one-dimensional case, let $\theta$ be position, and $\xi$ be momentum. These are both maximally accessible by the Heisenberg inequality. Let $G$ be the translation group acting on $\theta$. Then (i) is satisfied, and (ii) can be shown to be satisfied. The variable $\phi$ may be the vector $(\theta, \xi )$, and based on this vector, the transformation $k$ may be defined in terms of a Fourier transform: $(\theta, \xi)\mapsto (F[\theta], F^{-1}[\xi])$, where $F$ is the relevant Fourier transform.

The same theorem can also be used as a basis for a quantum description of a decision situation, if we accept the model descibed in the previous Section. We then must accept the technical assumptions (i)-(iii), but we are helped by the fact that there already exists a Quantum Decision Theory, founded by Yukalov and Sornette (2014) and related papers.

The theory is general. Operators corresponding to accessible conceptual variables that are not maximally accessible, can be defined by the spectral theorem (Hall, 2013), by taking as a point of departure that if $\theta$ is any accessible variable, then there exists a maximally accessible $\lambda$ and a function $f$ such that $\theta=f(\lambda)$.

The textbook treatment of quantum mechanics concentrates usually on the case where the space $\Omega_\theta$ is finite. From the theory of Helland, 2021 it is deduced: The operator $A^\theta$ has a finite spectrum. The set of  eigenvalues of $A^\theta$ coincides with the set of possible values for $\theta$. The eigenspaces of operators corresponding to a set of conceptual variables $\theta$ are in one-to-one correspondence with a question `What will $\theta$ be if I measure it?' together with a sharp answer `$\theta=u$'. Finally, $\theta$ is maximally accessible if and only if all eigenspaces of $A^\theta$ are one-dimensional.

The assumptions of Theorem 1 can tentatively be made more explicit in this finite case. To be precise, the assumption (ii) of Theorem 1 seem to follow in this case (Burnside, 1955) if one can prove that the actual representation $V$ satisfies the following: Every irreducible representation of $G$ occurs as a subrepresentation of $V^{\otimes n}$ for sufficiently large $n$. In general, we need the existence of a group $G$ and the corresponding representation $V=\{U(\cdot)\}$ satisfying the assumptions (ii) and (iii).  When the dimension of $\theta$ is $2p$ for some $p$, we can use a group  corresponding to $p$ independent qubits, and the necessary results are known (H\"{o}hn, 2017; H\"{o}hn and Wever, 2017).

I will not go into these technicalities here, but conclude: The textbook treatment of quantum mechanics can be simplified if one starts with the notion of conceptual variables. The operator corresponding to an accessible variable can be explicitly defined (Helland, 2021). The eigenvalues of these operators have a simple interpretation. Very many unit vectors in the Hilbert space also have a simple interpretation: They can be interpreted in terms of a question-and-answer pair. To characterize situations where this interpretation is valid for \emph{all} unit vector in the Hilbert space, is a difficult problem, stated as a question to the quantum community in Helland (2019).

In order to do calculations, we have to express the above in a physical language. Let $\theta$ be a maximally accessible conceptual variable taking a finite number of values $\{u_i ; i=1,...,d\}$. Then the corresponding Hilbert space is $d$-dimensional, and can be taken as $\mathbb{C}^d$, the vector space of $d$-dimensional complex vectors. Following Dirac, the (coloumn) vectors in this space are denoted as $|\psi\rangle$, called ket vectors. The corresponding complex conjugate row vectors are denoted by $\langle\psi |$, and are called bra vectors. The scalar product can then be defined as $\langle\psi_1 |\psi_2\rangle$.

Let $A^\theta$ be the operator (matrix) corresponding to $\theta$, denote the eigenvectors of $A^\theta$ as $\{|\psi_i\rangle ; i=1,...,d\}$, and the corresponding eigenvalues by $\{u_i ;i=1,...,d\}$. Then $\{|\psi_i \rangle\}$ is an orthonormal basis of the Hilbert space $\mathbb{C}^d$, and each state vector $|\psi_i\rangle$ is associated with the event $\theta=u_i$.

The state vector $|\psi_i\rangle$ may equally well be expressed by the projection operator $P_i =|\psi_i\rangle\langle\psi_i |$. The $P_i$'s are orthogonal, and satisfy $\sum_i P_i =I$, the identity.

Let now $\eta$ be another accessible conceptual variable, not necessarily maximal. Then $\eta$ is associated with an operator $A^\eta$, which may have degenerate eigenvalues $v_i$ $(i=1,...,r)$, where $r\le d$. And again the eigenspace of $A^\eta$ corresponding to the eigenvalue $v_i$ can be associated with the statement $\eta=v_i$. Let $P_i$ be the projection upon this eigenspace. Then again the different $P_i$ are orthogonal, and $\sum P_i =I$. Each $P_i$ can, in several ways, be written as a sum of one-dimensional projections on orthogonal eigenvectors.

In the quantum mechanical literature, a ket vector $|\psi\rangle$ is in general called a state vector, and it characterizes the state of a physical system. In my notation, this usually corresponds to a statement of the type $\theta =u$, where $\theta$ is maximally accessible. The observer $A$, being in this state relative to the physical system, is then certain of the value of $\theta$. Note that $\theta$ can be a vector, whose components are independent, accessible conseptual variables.

In other cases, $A$ may be uncertain of the value of $\theta$. A special case of this again, is that he for each $i$ has a probability $p_i$ for each value $\theta=u_i$. This is also by physicists considered to be a state, a mixed state $\rho$. Explicitly
\begin{equation}
\rho = \sum_i p_i P_i ,
\label{11}
\end{equation}
where $P_i$ is the projector onto the eigenspace of $A^\theta$ corresponding to $\theta = u_i$.

Several possibilities exist for the probabilities $p_i$. Depending upon the situation, they can be Bayesian priors or posteriors, or they can also have a frequentist interpretation, being derived from a fiducial / confidence distribution.
Other possibilities may exist.

It can be verified: $\rho$ is positive definite and has trace 1. In general, a density operator, characterizing a mixed state is defined as any operator with these two properties.

To complete the foundation of quantum mechanics, we need a device to calculate probabilities. This is given by Born's formula, the dervation of which is derived by several researchers from different assumptions. One such derivation is in Chapter 5 of Helland (2021).

\section{Discussion}
\label{Sec11}

As a general point, it is said in Helland (2021) that the scientific development, both within statistical inference and within quantum physics would perhaps have been slower if a connection between the two cultures had been taken into account from the beginning. Universality and creativity can in some sense be seen as complementary properties.

The discussion in this article should be contrasted to the long history that lies behind the development of quantum machanics, see Gilder (2008). Many great minds have arrived at different interpretations of the theory. The present discussion has been connected to the process of making decisions. Real decisions are much more complicated than what has been outlined above. In particular, this applies to decisions made by scientists, both physicists and statisticians. For these groups, as for other groups, decisions have a strong cultural component.

Political leaders around the world also make decisions. Some of these decisions are real cruel and reprehensible, and should be opposed by all thinking people. Unfortunately, in some connections the concept of truth seems to be seen as relative, in contrast to all aims of science. A common goal for scientists and politicians may be formulated by good leaders. These goals should include basic humanitarian values and considerations for future generations.

To take a much more modest problem, what can statisticians learn by considering the discussion above? Look first at the PLS modeling case. Here it seems to be clear that the PLS algorithm has a role to play in statistical predictions from regression models with near collinear data. More specifically (Cook and Forzani, 2019) the algorithm seems to give good predictions in cases where, in contrast to the sparsity case, several  combinations of $x$-variables play together in their effect on the response $y$.

Concerning a possible link between quantum theory and macroscopic activities like statistical inference, it may perhaps be too early to say much about what we can learn. Again, some tentative discussion is given in Helland (2021). As a strong simplification, it may perhaps be of some value for statisticians to look at the notion of taking decisions in the way that is described above. A problem area where this may be investigated, is in connection to the multiplicity issue in relation to the interpretation of p-values, an area where there has been much controversy recently (Wasserstein et al., 2019; Benjami et al., 2021).

In any case, it may be of some value to contrast the use of mathematical modeling in different areas of science.

\section*{References}

Benjamini, Y., De Veaux, R.D., Efron, B., Evans, S., Glickman, M., Graudbard, B.I., He, X., Meng, X.-L., Reid, N. Stigler, S.M., Vardeman, S.B., Wikle, C.K., Wright, T., Young, L.D. \& Kafadar, K. (2021). The ASA president's task force statement on statistical significance and replicability. \textit{Annals of Applied Statistics} 15 (3), 1084-1085.

Bhadra, A., Datta, J., Polson, N., Sokolov. V. and Xu, J. (2021). Merging two cultures: Deep and statistical learning. arXiv: 2110.11561 [stat.ME].

Breiman, L. (2001). Statistical modeling: The two cultures. \textit{Statistical Science} 16, 199-231.

Burnside, W. (1955). \textit{Theory of groups of finite order.} Dover Publications, Inc., New York.

Cook, R.D. (2018). \textit{An Introduction to Envelopes.} Wiley, Hobroken, N.J.

Cook, R.D., Li, B. and Chiaromonte, F. (2010). Envelope models for parsimonious and efficient multivariate linear regression. \textit{Statistica Sinica} 20, 927-1010.

Cook, R.D., Helland, I.S. and Su, Z. (2013). Envelopes and partial least squares regression. \textit{Journal of the Royal Statistical Society} B. 75, 851-877.

Cook, R.D. and Forzani. L. (2019). Partial least squares prediction in high-dimensional regression.  \textit{Annals of Statistics} 47 (2), 884-908.

Cox, D.R. and Wermuth, N. (2004). Causality: A statistical view. \textit{International Statistical Review} 72 (3) 285-305.

Efron, B. and Hastie, T. (2021). \textit{Computer Age Statistical Inference. Algorithms, Evidence, and Data Science.} Student Edition. Cambridge University Press, Cambridge.

Fisher, R.A. (1973). \textit{Statistical methods and scientific inference.} Hafner Press, New York.

Fraser, D.A.S. (1964). On the definition of fiducial probability. \textit{Bulletin of the International Statistical Institute} 40, 842-856.

Gilder, L. (2008). \textit{The Age of Entanglement. When Quantum Physics was Reborn.} Vintage Books, New York.

Hall, B.C. (2013). \textit{Quantum Theory for Mathematician.} Springer, New York.

Helland, I.S. (1990). Partial least squares and statistical models. \textit{Scandinavian Journal of Statistics} 17, 97-114.

Helland, I.S. (1992). Maximum likelihood regression on relevant components. \textit{Journal of the Royal Statistical Society} B. 54 (2), 637-647.

Helland, I.S: (2010). \textit{Steps Towards a Unified Basis for Scientific Models and Methods.} World Scientific, Singapore.

Helland, I.S. (2019). When is a set of questions to nature together with sharp answers to those questions in one-to-one correspondence with a set of quantum states? arXiv: 1909.08834 [quant-ph].

Helland, I.S. (2021). \textit{Epistemic Processes. A Basis for Statistics and Quantum Theory.} Springer, Cham, Switzerland.

Helland, I.S, (2022a). On reconstructing parts of quantum theory from two related maximal conceptual variables.  \textit{International Journal of Theoretical Physics} 61, 69.

Helland, I.S. (2022b). The Bell experiment and the limitations of actors. \textit{Foundations of Physics} 52, 55..

Helland, I.S. and Alm\o y, T. (1994). Comparison of prediction methods when only a few components are relevant. \textit{Journal of the American Statistical Association} 89 (426), 583-591.

Helland, I.S., S\ae b\o , S. and Tjelmeland, H. (2012). Near optimal prediction from relevant components. \textit{Scandinavian Journal of Statistics} 39, 695-713.

H\"{o}hn, P.A. (2017). Quantum theory from rules on information acquisition. \textit{Entropy} 19 (3), 98.

H\"{o}hn, P.A. and Wever, C.S.P. (2017). Quantum theory from questions. \textit{Physics Reviews} A, 95, 012102.

J\"{o}reskog, K.G. and Wold, H. (1982). \textit{Systems under indirect observation. Causality-structure-prediction.} Part I-II. North-Holland, Amsterdam.

Khan, M.E. and Rue, H. (2021) The Bayesian learning rule. arXiv: 2107.04562 [stat.ML].

Khare, K., Pal, S. and Su, Z. (2017) A Bayesian approach to envelope models. \textit{Annals of Statistics} 45 (1), 196-222.

Martens, H. (1987). Why users like PLS regression.  In: Martens, M. (ed.). \textit{Data approximation by PLS methods.} Report no. 800, Norwegian Computing Centre.

Mehmood, T. and Ahmed, B. (2015). The diversity in the applications of partial least squares: an overview. \textit{Journal of Chemometrics} 30 (1), 4-17.

N\ae s, T. and Helland, I.S. (1993). Relevant components in regression. \textit{Scandinavian Journal of Statistics} 20, 239-250.

Pearl, J. (2009). \textit{Causality. Models, Reasoning, and Inference.} 2. edition. Cambridge University Press, Cambridge.

Pothos, E.M. and Busemeyer, J.R. (2022). Quantum cognition. \textit{Annual Reviews of Psychology} 73, 749-778.

Schweder, T. (1987). Canonical regression versus PLS. In: Martens, M. (ed.). \textit{Data approximation by PLS methods.} Report no. 800, Norwegian Computing Centre.

Schweder, T and Hjort, N.L. (2016). \textit{Confidence, Likelihood, Probability, Statisical Inference with Confidence Distributions.} Cambridge University Press, Cambridge.

Taraldsen, G. and Lindqvist, B.H. (2021). Fiducial inference and decision theory. arXiv: 2112.07060 [stat.ME].

Yukalov, V.I. and Sornette, D. (2014). How brains make decisions. \textit{Springer Procedings in Physics} 150, 37-53.

Yukalov, V.I., Yukalova, E.P. and Sornette, D. (2018). Information processing by networks of quantum decision makers. \textit{Physics} A 482, 747-766.

Yukalov, V.I., Yukalova, E.P. and Sornette, D. (2022). Role of collective information in networks of quantum operating systems. arXiv: 2201.11008 [physics.soc-ph].

Wasserstein, R.L., Schirm, A.L. \& Lazar, N.A. (2019). Moving to a world beyond `$p<0.05$'. Editorial. \textit{American Statistician} 73: sup1: 1-19.

Zwirn, H. (2016). The measurement problem: Decoherence and convivial solipsism. \textit{Foundations of Physics} 46, 635-667.

\end{document}